\hfuzz=6pt
\font\titlefont=cmbx10 scaled\magstep1

\magnification=1200
\line{}
\vskip 1.5cm
\centerline{\titlefont MORE ON THE q-OSCILLATOR ALGEBRA}
\smallskip
\centerline{\titlefont AND q-ORTHOGONAL POLYNOMIALS}
\vskip 2cm
\centerline{\bf Roberto Floreanini}
\smallskip
\centerline{Istituto Nazionale di Fisica Nucleare, Sezione di Trieste}
\centerline {Dipartimento di Fisica Teorica,
Universit\`a di Trieste}
\centerline{Strada Costiera 11, 34014 Trieste, Italy}
\vskip 1cm
\centerline{\phantom{$^{(*)}$}{\bf Jean LeTourneux}$^{(*)}$}
\medskip
\centerline{and}
\medskip
\centerline{\phantom{$^{(*)}$}{\bf Luc Vinet}\footnote{$^{(*)}$}{Supported 
in part by the National Sciences and Engineering
Research Council \hbox{(NSERC)} of Canada and the Fonds FCAR of Qu\'ebec.}}
\smallskip
\centerline{Centre de Recherches Math\'ematiques}
\centerline{Universit\'e de Montr\'eal}
\centerline{Montr\'eal, Canada H3C 3J7}
\vskip 1.5cm
\centerline{\bf Abstract}
\smallskip\midinsert\narrower\narrower\noindent  
Properties of certain $q$-orthogonal polynomials are connected to the 
$q$-oscillator algebra. The Wall and $q$-Laguerre polynomials are shown to 
arise as matrix elements of $q$-exponentials of the generators in a 
representation of this algebra. A realization is presented where the
continuous $q$-Hermite polynomials form a basis of the representation space.
Various identities are interpreted within this model. In particular,
the connection formula between the continuous big $q$-Hermite polynomials
and the continuous $q$-Hermite polynomials is thus obtained, and two
generating functions for these last polynomials are algebraically derived.
\endinsert

\vfill\eject

The $q$-oscillator algebra is generated by the elements $A_+$, $A_-$
and $K = q^{-N/2}$ subjected to the relations
$$\eqalign{
&A_-\,A_+ - {1\over q} A_+\, A_- = 1\ ,\cr
&KA_+ = q^{-1/2} A_+ K\ ,\cr
&KA_- = q^{1/2} A_- K\ .} \eqno(1)$$
It represents, clearly, a $q$-deformation of the oscillator algebra which is
retrieved in the limit $q\rightarrow 1$.  The algebra (1) has been found to have
a number of applications.  It was shown in particular to provide the
algebraic interpretation of various $q$-special functions.  
(See for instance Refs.[1-11].)  The purpose
of the present note is to present additional results on this topic.

Most special functions have $q$-analogs;$^{12,13}$ these arise in various
connections and, in particular, in the description of systems with quantum
group symmetries.  The algebraic interpretation of many $q$-special functions
has been shown to proceed
in analogy with the Lie theory treatment of their classical, $q\rightarrow 1$,
counterparts. (See for instance Refs.[1, 9, 11].)

One considers $q$-exponentials of the generators of a $q$-algebra and observes
that their matrix elements in representation spaces are expressible in terms of
$q$-special functions. One then uses models to derive properties of these
functions through symmetry techniques. This approach has 
proven very fruitful and is still actively pursued.  
It is the one that we shall adopt here.

The outline of this letter is the following. We shall first introduce the
representation of (1) that will be used, and give a realization of this
representation where the continuous $q$-Hermite polynomials play the role of
the basis vectors.  Matrix elements of $q$-exponentials of $A_+$ and $A_-$ will then be evaluated and
seen to involve, in some cases, Wall and $q$-Laguerre polynomials. Finally,
these results will be used in conjunction with the model we mentioned to
obtain formulas involving the continuous $q$-Hermite polynomials.  
A connection formula and two generating
function identities will thus be algebraically derived.

We shall consider the following representation of the $q$-oscillator algebra. 
We shall denote by $\xi_n$, $n = 0,1,2,\dots$ the basis vectors, and take the
generators $A_+$, $A_-$ and $K$ to act according to
$$\eqalign{
&A_+ \xi_n
	= -q^{-(n+1)/2}\, \xi_{n+1}\ ,\cr
&A_- \xi_n
	= q^{n/2+1} \left( {1-q^{-n}\over 1-q} \right) \xi_{n-1}\ ,\cr
&K \xi_n
	= q^{-n/2}\, \xi_n\ .}\eqno(2)$$
It is immediate to check that these definitions are compatible with the
commutation relations (1). A first connection with $q$-polynomials is
made by observing that there exists a realization of this representation where
the continuous $q$-Hermite 
polynomials $H_n(x|q)$ appear as basis vectors.$^{14}$  These polynomials
are defined as follows:$^{13}$
$$H_n(x|q) = \sum_{k=0}^n \left[{n\atop k}\right]_q e^{i(n-2k)\theta}\ , \quad
x = \cos \theta\ , \eqno(3)$$
with the $q$-binomial coefficients given by
$$\left[{n\atop k}\right]_q =\left[{n\atop n-k}\right]_q =
{(q;q)_n \over(q;q)_k\, (q;q)_{n-k}}\ . \eqno(4)$$
We use standard notation$^{12,13}$ where $(a;q)_\alpha$ stands for
$$(a;q)_\alpha={(a;q)_\infty \over(aq^\alpha; q)_\infty}\ ; 
\quad (a;q)_\infty =\prod_{k=0}^\infty (1-aq^k), \quad |q| < 1\ . \eqno(5)$$
The classical Hermite polynomials $H_n(x)$, are obtained as follows from the
continuous $q$-Hermite polynomials $H_n(x|q)$ in the limit $q\rightarrow 1$:
$$
\lim_{q \rightarrow 1^-} \left({1-q\over 2}\right)^{-n/2} H_n\big( x
\sqrt{(1-q)/2} \big| q \big) = H_n (x)\ . \eqno(6)$$
Let $T_z$ be the $q$-shift operator:
$$
T_z\, f(z) = f(qz)\ . \eqno(7)
$$
It can be verified$^{15}$ that the representation $(2)$ is realized by
setting$^{14}$
$$
\xi_n(x) = H_n(x|q), \qquad n=0,1,2, \dots\ , \eqno(8)
$$
and taking $A_+$, $A_-$ and $K$ to be the following operators acting on
functions of $x = (z + z^{-1})/2$, $z = e^{i\theta}$:
$$
\eqalignno{
&A_+= q^{-1/2} {1\over z-z^{-1}} \left({1\over z^2}\, 
T_z^{1/2} - z^2T_z^{-1/2}\right)\ ,&(9a)\cr
&A_-= {q\over 1-q} \ {1\over z-z^{-1}} (T_z^{1/2} - T_z^{-1/2})\ , 
&(9b)\cr
&K= {1\over z-z^{-1}} \left( -{1\over z}\, 
T_z^{1/2} + z\, T_z^{-1/2} \right)\ . &(9c)}
$$
The operator $\tau = (z - z^{-1})^{-1}(T_z^{1/2} - T_z^{-1/2})$ in $A_-$ is
referred to as the divided difference operator.

In order to mimic the exponential mapping from Lie algebras to Lie groups, we
shall need $q$-analogs of the function $e^x$. Introduce the $q$-exponentials
$$
E_q^{(\mu)} (x)
	= \sum_{n=0}^\infty {q^{\mu n^2}\over (q;q)_n} x^n, \qquad \mu \in {\bf R}
\ .\eqno(10)
$$
In the limit $q \rightarrow 1$, once $x$ has been rescaled by $(1-q)$, these 
functions all tend to the ordinary exponential: $\lim_{q \rightarrow 1} 
E_q^{(\mu)}\bigl[ (1-q) x \bigr] = e^x$.  
For some specific values of $\mu$, they correspond to standard $q$-exponentials.
Indeed, for $\mu = 0$ and $\mu = 1/2$ one has$^{12}$
$$
\eqalignno{
&E_q^{(0)} (x)
	=e_q(x) = {1\over (x;q)_\infty}\ , &(11a)\cr
&E_q^{(1/2)} (x)
	= E_q(q^{-1/2} x) = (-q^{-1/2}x; q)_\infty\ . &(11b)}
$$
Note that $e_q(\lambda x)$ and $E_q(-q\lambda x)$ are, respectively,
eigenfunctions with eigenvalue $\lambda$ of the $q$-derivative operators 
$D_z^+ =z^{-1}(1-T_z)$ and $D_z^- = z^{-1}(1 - T_z^{-1})$.  
Not so well known is the $q$-exponential
$$
{\cal E}_q(x;a,b) =\sum_{n=0}^\infty {q^{n^2/4}\over (q;q)_n} 
\bigl( aq^{(1-n)/2} e^{i\theta};q \bigr)_n  \,
\bigl( aq^{(1-n)/2} e^{-i\theta}; q \bigr)\,_n\, b^n\ ,\quad  x = \cos\theta\ ,
\eqno(12)
$$
introduced in Ref.[16]. It enjoys the property of being an eigenfunction of the
divided difference operator $\tau$ with eigenvalue $ab \,q^{-1/4}$, and, in the
limit $q \rightarrow 1$, 
${\cal E}_q \bigl(x;a,(1-q)b \bigr) \rightarrow \exp \bigl[ (1 + a^2 -2ax) b 
\bigr]$. We observe that
$$
E_q^{(1/4)} (x) = {\cal E}_q(-; 0, x)\ .\eqno(13)
$$
The operators
$$
U^{(\mu, \nu)} (\alpha, \beta) =
	E_q^{(\mu)} \bigl( (1-q)\alpha\, A_+ \bigr)\  
E_q^{(\nu)} \left( {\beta\over q} (1-q)\,A_- \right)\ ,\eqno(14)
$$
in the completion of the $q$-oscillator algebra, are central in our analysis. 
In the limit $q\rightarrow 1$, 
they go into the Lie group element $e^{\alpha A_+}e^{\beta A_-}$.  
Their matrix elements in the representation space spanned by
the vectors $\xi_n$ are defined by
$$
U^{(\mu, \nu)}(\alpha, \beta)\,  \xi_n
	= \sum_{m=0}^\infty U_{m,n}^{(\mu,\nu)} (\alpha, \beta)\, \xi_m\ ,\eqno(15)
$$
and, when evaluated, are found to involve $q$-special functions.

Explicitly, one obtains
$$
\eqalignno{
&U_{m,n}^{(\mu,\nu)}(\alpha,\beta)
	= (-\beta)^{n-m}q^{(n-m)[(\nu +1/4)(n-m)-n/2-1/4]}
\left[{n\atop m}\right]_q  &(16a)\cr
	&\hskip 4cm \times {\cal P}_m^{(\mu,\nu)} \bigl( - (1-q)\alpha \beta;
q^{n-m} | q\bigr), \qquad {\rm if } \ m \leq n\ ,\cr
&U_{m,n}^{(\mu,\nu)}(\alpha,\beta)
	={[-(1-q)\alpha]^{m-n}\over(q;q)_{m-n}}
q^{(n-m)[(\mu-1/4)(m-n)-n/2-1/4]} &(16b)\cr
	&\hskip 4cm\times {\cal P}_n^{(\nu,\mu)} \bigl( - (1-q)\alpha \beta;
q^{m-n} | q\bigr), \qquad {\rm if } \ m \geq n\ ,}
$$
where ${\cal P}_n^{(\mu,\nu)} (x;q^\gamma | q)$ are the polynomials given by
$$
{\cal P}_n^{(\mu,\nu)}(x;q^\gamma |q) =
	\sum_{k=0}^n {q^{k^2(\mu+\nu) + 2\nu \gamma k}(q^{-n};q)_k\over
(q;q)_k\, (q^{\gamma +1};q)_k}\ x^k\ . \eqno(17)
$$
Note that in passing from expression $(16a)$ to $(16b)$ for
$U_{m,n}^{(\mu,\nu)} (\alpha,\beta)$ or vice-versa, $m$ and $n$ as well as
$\mu$ and $\nu$, are exchanged in the polynomials ${\cal P}_n^{(\mu,\nu)}$.

The connection with standard $q$-polynomials is observed for particular values
of $\mu$ and $\nu$.  The little $q$-Laguerre or Wall polynomials$^{13}$
$$
p_n(x; a|q) = {}_2 \phi_1 \left({q^{-n},\atop\ }{0\atop aq}\bigg|\, q; 
q x \right)
\ ,\eqno(18)
$$
are, for instance, seen to occur for $\mu = \nu = 0$.  Indeed, 
$$
{\cal P}_n^{(0,0)} 
\bigl( -(1-q)\alpha \beta; q^{m-n} | q \bigr) = p_n \bigl( 
(1- 1/q)\alpha\beta; q^{m-n} | q\bigr)\ .\eqno(19)
$$
Similarly, the $q$-Laguerre polynomials$^{12,13}$
$$
L_n^{(\rho)}(x;q) ={(q^{\rho + 1};q)_n\over (q;q)_n}\ {}_1\phi_1
\left({q^{-n}\atop q^{\rho+1}}
\bigg|\, q ; -xq^{n+\rho+1} \right)\ ,\eqno(20)
$$
are found to arise when $\mu = \nu = 1/4$.  In this case we find for example
$$
{\cal P}_n^{(1/4, 1/4)} \bigl(x; q^{m-n} | q \bigr) 
={(q;q)_n\over(q^{m-n+1};q)_n} L_n^{(m-n)} \bigl(x; q^{-(m+n+1)/2};q \bigr)\ .
\eqno(21)
$$
The $q$-hypergeometric series ${}_r \phi_s$ that we are using are defined 
by$^{12,13}$
$$
\eqalign{
{}_r \phi_s &\left[{a_1,\atop b_1,} {a_2,\atop b_2,} {\dots\atop\dots}
{a_r\atop b_s} \bigg|\, q; z\right]\cr
&=\sum_{n=0}^\infty {(a_1;q) \dots (a_r;q)_n\over (q;q)_n(b_1,q)_n \dots
(b_s;q)_n} \bigl[ (-1)^n q^{n(n-1)/2} \bigr]^{1+s-r}\ z^n\ .}\eqno(22)
$$

We shall now return to the model (8), (9) and make use of these
matrix elements to derive properties of the continuous $q$-Hermite polynomials.
We shall first produce a formula giving the expansion of the continuous big
$q$-Hermite polynomials $H_n(x;a|q)$ in terms of 
continuous $q$-Hermite polynomials $H_n(x|q)$; that
is, we shall derive from symmetry considerations an identity of the form
$$
H_n(x;a|q) = \sum_{k=0}^\infty C_{k,n} H_k(x|q)\ .\eqno(23)
$$
The continuous big $q$-Hermite polynomials are defined by$^{13}$
$$
H_n(x;a|q) =
	e^{in\theta}\ {}_2\phi_0
\left({q^{-n},\ ae^{i\theta}\atop -}\bigg|\,q; q^n e^{-2i\theta}
\right)\ , \qquad x = \cos \theta\ . \eqno(24)
$$
An algebraic interpretation of these polynomials is found in Ref.[17].  It is
easy to see from the definition (3) that the continuous $q$-Hermite
polynomials $H_n(x|q)$ are the $a \rightarrow 0$ limits of the big ones, 
{\it i.e.} $H_n(x|q) = H_n(x;0|q)$. Since the $H_n(x;a|q)$ and {\it a fortiori}
the $H_n(x|q)$ are particular cases of
Askey-Wilson polynomials, the coefficients $C_{k,n}$ in (23) can
evidently be obtained by specializing the general formula for the connection
coefficients of these polynomials.$^{12}$ 
Our purpose here is to show that they can be given an algebraic interpretation.

We adopt the realization (8), (9) of the representation (2) for the 
$q$-oscillator algebra, and
consider within this framework the action of 
${\cal E}_q(-;0, (\beta(1-q)/q) A_-)= U^{(0,1/4)}(0,\beta)$ 
on continuous $q$-Hermite polynomials.  Since the computation of the
matrix elements of $U^{(\mu, \nu)}(\alpha, \beta)$ is model independent, we
have, on the one hand, from (15) and $(16a)$:
$$
\eqalign{
{\cal E}_q & \Big(-;0,(\beta(1-q)/q)\, A_- \Big)\, H_n(x|q) \cr
&\hskip 2cm =\sum_{k=0}^n (-1)^k a^k q^{k(k-1)/2} \left[{n\atop k}\right]_q 
H_{n-k}(x|q)\ ,}\eqno(25) 
$$
with
$$ a = q^{-n/2 + 1/4} \beta\ .\eqno(26)$$
On the other hand, it turns out that the action of ${\cal E}_q(-;0,
(\beta(1-q)/q) A_-)$ on $H_n(x|q)$ can be resummed, using the explicit
expression (3) and the $q$-binomial theorem,
$$
\sum_{n=0}^\infty {(\alpha ; q)_n\over (q;q)_n} z^n =
{(\alpha z;q)_\infty\over (z;q)_\infty}\ . \eqno(27)
$$
One arrives at
$$
U^{(0,1/4)}(0,\beta)\, H_n(x|q) =
(a/z;q)_n z^n\, {}_2 \phi_1 \left({q^{-n},\atop\ }{0\atop q^{1-n} z/a}
\bigg|\, q ; q/az\right)\ .\eqno(28)
$$
One then uses the transformation formula$^{13}$
$$
\eqalign{
{}_2\phi_1 \left({q^{-n},\atop\ }{0\atop c}\bigg|\, q ; z\right)=&
(-1)^n q^{-n(n+1)/2} {z^n\over (c;q)_n}\cr
&\times {}_2\phi_0\left({q^{-n},\ q^{1-n}/c\atop -} \bigg|\, q ; 
q^{2n}c/z\right)\ ,}\eqno(29) 
$$
and the definition (26) to show that:
$$
{\cal E}_q(-;0,\beta \tau)\, H_n(x|q) = H_n(x;a|q)\ ,\eqno(30)
$$
where $a$ is still given by (26). Putting (25) and (30) together
finally yields the expansion formula:
$$
H_n(x;a|q) =
	\sum_{k=0}^n(-1)^k a^k q^{k(k-1)/2} \left[{n\atop k}\right]_q\
H_{n-k}(x|q)\ .\eqno(31)
$$
In the limit $q \rightarrow 1$, this relation tends to the following 
identity between classical Hermite polynomials:
$$
H_n(x-a) = \sum_{k=0}^n(-1)^{n-k} (2a)^{n-k} \left({n\atop k}\right)\ 
H_k(x)\ .\eqno(32)
$$
This can be verified using (6) and noting that
$$
\lim_{q \rightarrow 1^-} \left({1-q\over 2}\right)^{n/2}\ H_n \left(
x \sqrt{(1-q)/2}; a \sqrt{2(1-q)} \bigg|\, q \right) = H_n (x-a)\ .\eqno(33)
$$
As a second example of application, we shall constructively derive two
generating function identities.  Consider the action of the operators
$U^{(\mu,\nu)} \bigl( \alpha/(1-q), 0 \bigr)$ on $\xi_0(x) = 1$.  We have from
(15) and $(16b)$
$$
U^{(\mu,0)}\bigl( \alpha/(1-q), 0 \bigr) \cdot 1 = 
\sum_{m=0}^\infty {(-q^{-1/4} \alpha)^m\over (q;q)_m} q^{(\mu - 1/4)m^2}\ 
H_m(x|q)\ . \eqno(34)
$$
On the one hand, this relation tells how the function $\xi_0(x) = 1$ transform
under the action of what would be in the limit $q \rightarrow 1$, a 
transformation group element.  
On the other hand, if expressions in closed form can be found
for $U^{(\mu,0)} \bigl( \alpha / (1-q), 0 \bigr) \cdot 1$, 
these would be generating functions for
the continuous $q$-Hermite polynomials.

It is readily seen that two such generating functions can indeed derived from
(34) for $\mu = 1/4$ and $\mu = 3/4$.  In the first case one has, after
having set $\mu = 1/4$ and inserted the explicit expansion of $H_m(x|q)$ in
(34):
$$
{\cal E}_q(-;0, \alpha A_+) \cdot 1 =
	\sum_{m,k=0}^\infty {(-\alpha q^{-1/4})^m\over(q;q)_{m-k}(q;q)_k}\
z^{m-2k}\ .\eqno(35)
$$
When using $\ell = m-k$ instead of $m$ as summation index, the two sums are
seen to split and one finds that
$$
{\cal E}_q(-;0,\alpha A_+) \cdot 1 = e_q(-q^{-1/4} \alpha z)\ 
e_q \left( -q^{-1/4} \alpha/z \right)\ , \eqno(36)
$$
where $e_q(x)$ is the $q$-exponential defined in $(11a)$. One now sets 
$t= -q^{-1/4}\alpha$ to see that (34) entails in the case $\mu = 1/4$, the
generating function identity:
$$
e_q(tz)\ e_q(t/z) =
	\sum_{n=0}^\infty {t^n\over (q;q)_n}\ H_n(x|q)\ . \eqno(37)
$$
The action of $U^{(3/4,0)}(\alpha /(1-q),0)$ on $\xi_0(x) = 1$ can similarly be
cast in closed form.  Using again (3) in (34), with $\mu = 3/4$, we
now have
$$
E_q^{3/4} (\alpha A_+) \cdot 1 =
	\sum_{k,m=0}^\infty {q^{m(m-1)/2}\over (q;q)_k(q;q)_{m-k}}\,
(-q^{1/4}\alpha)^m\ z^{m-2k}\ .\eqno(38)
$$
The two sums are reorganized by using $\ell = m-k$ instead of $m$ as summation
index. This allows one to perform the sum over $\ell$ thanks to the explicit
expansion for $(11b)$, and one thus arrives at
$$
E_q^{(3/4)} (\alpha A_+) \cdot 1 =(q^{1/4} \alpha z;q)_\infty\
{}_1\phi_1 \left({0\atop q^{1/4}\alpha z}\bigg|\,q; q^{1/4}\alpha/z\right)\ .
\eqno(39)
$$
We set $t = q^{1/4} \alpha$ and combine (38) and (39) to find
another generating function identity:
$$ 
(tz; q)_\infty\  {}_1\phi_1 \left({0\atop tz}\bigg|\,q; t/z\right)
=\sum_{k=0}^\infty {q^{n(n-1)/2}\over(q;q)_n}\, (-1)^n t^n\ H_n(x|q)\ ,
\eqno(40)
$$
the algebraic interpretation of which stems from (34).

The results presented here illustrate once more the usefulness of the algebraic
interpretation of $q$-special functions.  It is remarkable that relations
(31), (37) and (40), like many other $q$-special functions
identities, have their origin in 
the representation theory of the $q$-oscillator algebra.

\vskip 2cm

\centerline{\bf Acknowledgements}
One of us (L.V.) greatly benefitted from the support and hospitality extended
to him by the Sezione di Trieste of the INFN. 
	
\vskip 2cm

\centerline{\bf References}

\item{1.} Floreanini, R. and Vinet, L., $q$-Orthogonal polynomials
and the oscillator quantum group, Lett. Math. Phys. {\bf 22} (1991), 45--54.
\smallskip
\item{2.} Biedenharn, L.~C., The quantum group $SU(2)_q$ and a $q$-analogue
of the boson operators, J. Phys. A {\bf 22}, (1989), L873--L878.
\smallskip
\item{3.} MacFarlane, A.~J., On $q$-analogues of the quantum harmonic
oscillator and the quantum group $SU(2)_q$, J. Phys. A {\bf 22} (1989),
4581--4588.
\smallskip
\item{4.} Atakishiyev, N.~M. and Suslov, S.~K., Difference analogs of the
harmonic oscillator, Theor. Math. Phys. {\bf 85} (1990), 1055--l062; A
realization of the $q$-harmonic oscillator, ibid. {\bf 87} (1991), 442--444.
\smallskip
\item{5.} Kalnins, E.~G., Manocha, H.L. and Miller, W., Models of
$q$-algebra representations: I. Tensor products of special unitary
and oscillator algebras, J. Math. Phys. {\bf 33} (1992), 2365-2383. 
\smallskip
\item{6.} Kalnins, E.~G, Miller, W., and Mukherjee, S., Models of
$q$-algebra representations: Matrix elements of the $q$-oscillator algebra,
J. Math. Phys. {\bf 34} (1993), 5333--5356. 
\smallskip
\item{7.} Kalnins, E.~G. and Miller, W., Models of $q$-algebra
representations: $q$-integral transforms and ``addition theorems'', J. Math.
Phys. {\bf 35} (1994), 1951--1975.
\smallskip
\item{8.} Zhedanov, A.~S., Nonlinear shift of $q$-Bose operators and
$q$-coherent states, J. Phys. A {\bf 24} (1991), L1129--L1132; Weyl shift of
$q$-oscillator and $q$-polynomials, Theor. Math. Phys. {\bf 94} (1993),
219--224; $Q$-rotations and other $Q$-transformations as 
unitary nonlinear automorphisms of quantum algebras, J.
Math. Phys. {\bf 34} (1993), 2631--2647.
\smallskip
\item{9.} Floreanini, R. and Vinet, L., Quantum algebras and $q$-special
functions, Ann. Phys. {\bf 221} (1993) 53--70.
\smallskip
\item{10.} Floreanini, R. and Vinet, L., Automorphisms of the $q$-oscillator
algebra and basic orthogonal polynomials, Phys. Lett. A {\bf 180} (1993)
393--401.
\smallskip
\item{11.} Floreanini, R. and Vinet, L.,  ${\cal U}_q(sl(2))$ and
$q$-special functions, in {\it Lie Algebras, Cohomology and New Applications
to Quantum Mechanics}, Contemp. Math. {\bf 160} (1994), 85-100 .
\smallskip
\item{12.}  Gasper, G. and Rahman, M., {\it Basic Hypergeometric Series},
(Cambridge University Press, Cambridge, 1990).
\smallskip
\item{13.} Koekoek, R. and Swarttouw, R.~F., The Askey-scheme of
hypergeometric orthogonal polynomials and its $q$-analogue, Report 94--05,
Delft University of Technology (1994).
\smallskip
\item{14.} Floreanini, R. and Vinet, L., A model for the continuous
$q$-ultraspherical polynomials, CRM--2233, Universit\'e de Montr\'eal (1995).
\smallskip
\item{15.} Kalnins, E.~G. and Miller, W., Symmetry techniques for $q$-series:
Askey-Wilson polynomials, Rocky Mountain J. Math. {\bf 19} (1989), 223--230.
\smallskip
\item{16.} Ismail, M.~E.~H. and Zhang, R., Diagonalization of certain
integral operators, Adv. Math. {\bf 109} (1994), 1--33.
\smallskip
\item{17.} Floreanini, R., LeTourneux, J. and Vinet, L., An algebraic
interpretation of the continuous big $q$-Hermite polynomials, CRM-2246,
Universit\'e de Montr\'eal (1995).

\bye